\documentclass[11pt,twoside]{amsart}
\usepackage{tabularx}
\usepackage{geometry}
\usepackage{amssymb,amsmath,amsthm,latexsym}
\usepackage{mathrsfs}
\usepackage{url}
\usepackage{amsfonts}
\usepackage{graphicx}
\usepackage{hyperref}

\usepackage{graphicx}
\usepackage{amssymb,amsthm,amsmath}
\usepackage{mathtools}

\newtheorem{theorem}{Theorem}[section]
\newtheorem{corollary}{Corollary}[section]

\newtheorem{remark}{Remark}[section]

\usepackage{enumerate}

\usepackage{hyperref}
\usepackage{xcolor}

\hypersetup{colorlinks=true,linkcolor=blue,urlcolor=red}
\vspace{5cm}
  
\begin{document}

\title{Prime Divisors of 10's Friends: A Generalization of Prior Bounds} %%%%%%%%%%%%
\author[Sagar Mandal]{Sagar Mandal}
\address{Department of Mathematics and Statistics, Indian Institute of Technology Kanpur\\ Kalyanpur, Kanpur, Uttar Pradesh 208016, India}
\email{sagarmandal31415@gmail.com}

\maketitle
\let\thefootnote\relax
\let\thefootnote\relax
\maketitle

\begin{abstract}
10 is the smallest positive integer which is whether solitary or friendly is still an open question in mathematics. In this paper, we provide upper bounds for each of the prime divisors of a friend of 10. This paper is precisely a generalization of a recent paper \cite{4} in which necessary upper bounds for the 2nd, 3rd, and 4th smallest prime divisors of a friend of 10 have been proved. Further, we establish better upper bounds for the 3rd, and 4th smallest prime divisors of a friend of 10 than the bounds given in \cite{4}.
\end{abstract}

%\tableofcontents
\section{Introduction}
A positive integer $n$ other than $10$ is said to be a friend of 10 if $I(n)=I(10)=\frac{9}{5}$, where $I(n)$ is the abundancy index of $n$, which is defined as $I(n)=\sigma(n)/{n}$, where $\sigma(n)$ is the sum of positive divisors of $n$. The question "Is 10 a solitary number?" attracts many authors \cite{1,4,5,8,9}. In order to find a friend of 10, J. Ward \cite{9} proved that a friend $n$ of 10 is a square having at least six distinct prime divisors with $5$ being the smallest prime divisor, further, he described properties of some of the prime factors of $n$. In \cite{1}, the authors improved J. Ward's results by proving that a friend of 10 must have at least seven distinct prime divisors, additionally, they provided necessary properties of a friend of 10. In \cite{8}, the author proved that each friend of 10 has at least ten distinct prime divisors. In a recent paper \cite{4}, the authors proved the following upper bounds for the second, third, and fourth smallest prime divisors of friends of 10, \\
\begin{theorem}[\cite{4}, Theorem 1, 2, 3]\label{lem1}
If $n$ is a friend of 10 and if $q_2,q_3,q_4$ are the second, third, fourth smallest prime divisors of $n$ respectively. Then \\
\begin{align*}
7\leq q_2 < {\lceil\frac{7\omega(n)}{3}\rceil}\biggl\{\log\left({\lceil\frac{7\omega(n)}{3}\rceil}\right)+2\log\log\left({\lceil\frac{7\omega(n)}{3}\rceil}\right)\biggl\},\\
11\leq q_3 < {\lceil\frac{180\omega(n)}{41}\rceil}\biggl\{\log\left({\lceil\frac{180\omega(n)}{41}\rceil}\right)+2\log\log\left({\lceil\frac{180\omega(n)}{41}\rceil}\right)\biggl\},\\
13\leq q_4<  {\lceil\frac{390\omega(n)}{47}\rceil}\biggl\{\log\left({\lceil\frac{390\omega(n)}{47}\rceil}\right)+2\log\log\left({\lceil\frac{390\omega(n)}{47}\rceil}\right)\biggl\}.
\end{align*}
\end{theorem}
Where $\omega(n)$ is the number of distinct prime divisors of $n$ and $\lceil . \rceil$ is the ceiling function. In this paper, we modify these upper bounds and we give upper bounds for each of the prime divisors of a friend of 10 by proving the following theorem.
\begin{theorem}\label{thm 1.4}
   Let $n$ be a friend of 10 and $q_r$ ($r\geq 2$) be the $r$-th smallest prime divisor of $n$. Then we have \\
\begin{align*}
    q_r< L \{\log L+2\log(\log L)\},
\end{align*} 
where 
$$L=\lceil\frac{\mathcal{A}\omega(n)}{\mathcal{B}}\rceil,~~~~ \frac{\mathcal{A}}{\mathcal{B}} > \frac{1}{\frac{36}{25}\cdot {\displaystyle\prod_{4\leq i\leq r+1} (1-\frac{1}{p_i})-1}}~~~(\text{where}~p_i~\text{is the $i$-th prime number}),$$ 
$\mathcal{A}$, $\mathcal{B}$ are positive integers with $(\mathcal{A,\mathcal{B}})=1$
and $\frac{\mathcal{A}}{\mathcal{B}}\in \mathbb{Q}^{+}\setminus \mathbb{Z}^{+}$ such that $\mathcal{AB}(r-2)+2\mathcal{A}+\mathcal{B}>\mathcal{B}^2$.
\end{theorem}
\begin{remark}
    Since $\mathcal{A}>\mathcal{B}>1$, the condition $\mathcal{AB}(r-2)+2\mathcal{A}+\mathcal{B}>\mathcal{B}^2$ in the statement of Theorem \ref{thm 1.4} is redundant for $r\geq3$, and for $r=2$ it reduces to $2\mathcal{A}+\mathcal{B}>\mathcal{B}^2$.
\end{remark}
At the end of this note, we provide a proof of Theorem \ref{lem1} as a consequence of the above theorem. Further, as a corollary of the above theorem, we prove the following.
\begin{corollary}\label{coro1.4}
  Let $n$ be a friend of 10 and $q_3$, $q_4$ be the third and fourth smallest prime divisors of $n$ respectively. Then we have \\
\begin{align*}
11\leq q_3 < {\lceil\frac{427\omega(n)}{100}\rceil}\biggl\{\log\left({\lceil\frac{427\omega(n)}{100}\rceil}\right)+2\log\log\left({\lceil\frac{427\omega(n)}{100}\rceil}\right)\biggl\},\\
13\leq q_4<  {\lceil\frac{41\omega(n)}{5}\rceil}\biggl\{\log\left({\lceil\frac{41\omega(n)}{5}\rceil}\right)+2\log\log\left({\lceil\frac{41\omega(n)}{5}\rceil}\right)\biggl\}.
\end{align*}  
\end{corollary}
Observe that, $\mathcal{A}>\mathcal{B}>1$ and since a friend say $n$ of $10$ must have $\omega(n)>3$, thus $\lceil\frac{\mathcal{A}\omega(n)}{\mathcal{B}}\rceil>3$. Also, note that, $\log(\lceil\frac{\mathcal{A}\omega(n)}{\mathcal{B}}\rceil)+2\log\log(\lceil\frac{\mathcal{A}\omega(n)}{\mathcal{B}}\rceil)>0$, for  $\frac{\mathcal{A}\omega(n)}{\mathcal{B}}>3$, is an increasing function of $\frac{\mathcal{A}\omega(n)}{\mathcal{B}}$. Since $\lceil . \rceil$ is an increasing function, it follows that 
$$\lceil\frac{\mathcal{A}\omega(n)}{\mathcal{B}}\rceil \biggl\{\log(\lceil\frac{\mathcal{A}\omega(n)}{\mathcal{B}}\rceil)+2\log\log(\lceil\frac{\mathcal{A}\omega(n)}{\mathcal{B}}\rceil)\biggl\}$$
is an increasing function of $\frac{\mathcal{A}\omega(n)}{\mathcal{B}}$, for $\frac{\mathcal{A}\omega(n)}{\mathcal{B}}>3$. Since $\frac{427\omega(n)}{100}<\frac{180\omega(n)}{41}$ and $\frac{41\omega(n)}{5}<\frac{390\omega(n)}{47}$ we have
\begin{align*}
  q_3&< {\lceil\frac{427\omega(n)}{100}\rceil}\biggl\{\log\left({\lceil\frac{427\omega(n)}{100}\rceil}\right)+2\log\log\left({\lceil\frac{427\omega(n)}{100}\rceil}\right)\biggl\}  \\&< {\lceil\frac{180\omega(n)}{41}\rceil}\biggl\{\log\left({\lceil\frac{180\omega(n)}{41}\rceil}\right)+2\log\log\left({\lceil\frac{180\omega(n)}{41}\rceil}\right)\biggl\}  
\end{align*}
and 
\begin{align*}
  q_4&< {\lceil\frac{41\omega(n)}{5}\rceil}\biggl\{\log\left({\lceil\frac{41\omega(n)}{5}\rceil}\right)+2\log\log\left({\lceil\frac{41\omega(n)}{5}\rceil}\right)\biggl\}  \\&< {\lceil\frac{390\omega(n)}{47}\rceil}\biggl\{\log\left({\lceil\frac{390\omega(n)}{47}\rceil}\right)+2\log\log\left({\lceil\frac{390\omega(n)}{47}\rceil}\right)\biggl\}. 
\end{align*}
This shows that, Corollary \ref{coro1.4} provides better upper bounds for the prime divisors $q_3$ and $q_4$ of $n$ than the bounds given in Theorem \ref{lem1}.
\section{Properties of Abundancy Index}
The following are some of elementary properties \cite{3,10}  of the abundancy index.
\begin{enumerate}
    \item $I(F)$ is weakly multiplicative, that is, if $F$ and $G$ are two coprime positive integers then $I(FG)=I(F)I(G)$.
    \item\label{p2} If $\alpha,F$ are two positive integers and $\alpha>1$. Then $I(\alpha F)>I(F)$.
\item If $p_1$, $p_2$, $p_3,\dots ,p_n$ are $n$ distinct prime numbers and $\alpha_1$,$\alpha_2,\dots, \alpha_n$ are positive integers then
\begin{align*}
    I\biggl (\prod_{i=1}^{n}p_i^{\alpha_i}\biggl)=\prod_{i=1}^{n}\biggl(\sum_{j=0}^{\alpha_i}p_i^{-j}\biggl)=\prod_{i=1}^{n}\frac{p_i^{\alpha_i+1}-1}{p_i^{\alpha_i}(p_i-1)}.
\end{align*}
\item\label{p4} If $p_{1},\dots,p_{n}$ are distinct prime numbers and if $q_{1},\dots,q_{n}$ are distinct prime numbers such that  $p_{i}\leq q_{i}$ for all $1\leq i\leq n$ then for  positive integers $l_1,l_2,\dots,l_n$ we have
\begin{align*}
   I \biggl(\prod_{i=1}^{n}p_i^{l_i}\biggl)\geq I\biggl(\prod_{i=1}^{n}q_i^{l_i}\biggl).
\end{align*}
\item\label{p5}  If $F={\displaystyle\prod_{i=1}^{n}p_i^{\alpha_i}}$, then $I(F)<{\displaystyle\prod_{i=1}^{n}\frac{p_i}{p_i-1}}$.
\end{enumerate}
\vspace{5mm}
In order to prove our main theorem, we are required to use the following theorem.
\begin{theorem}[\cite{7}]\label{thm 3.4}
    If $p_n$ is $n\text{-th}$ prime number then,
    $$p_n<n(\log n+2\log\log n)$$
    for $n\geq 4$.
\end{theorem}
\section{Proof of Theorem \ref{thm 1.4}}
Let 
$$n=5^{2a_1}\cdot \prod_{2\leq i \leq \omega(n)}q^{2a_i}_i$$ 
be a friend of 10. At first, we shall prove that, $q_r$ must be strictly less than $p_{L}$ ($p_{L}$ is the $L$-th prime number) where
$$L=\lceil\frac{\mathcal{A}\omega(n)}{\mathcal{B}}\rceil,~~~\frac{\mathcal{A}}{\mathcal{B}} > \frac{1}{\frac{36}{25}\cdot {\displaystyle\prod_{4\leq i\leq r+1} (1-\frac{1}{p_i})-1}}~~~(\text{where}~p_i~\text{is the $i$-th prime number}),$$
$\mathcal{A}$, $\mathcal{B}$ are positive integers with $(\mathcal{A,\mathcal{B}})=1$
and $\frac{\mathcal{A}}{\mathcal{B}}\in \mathbb{Q}^{+}\setminus \mathbb{Z}^{+}$ such that $\mathcal{AB}(r-2)+2\mathcal{A}+\mathcal{B}>\mathcal{B}^2$.
Then by Theorem \ref{thm 3.4}, we shall have our desired result. If possible, suppose that, $q_r\geq p_{L}$. Since $q_2\geq 7=p_4,q_3\geq 11=p_5,\ldots, q_{r-1}\geq p_{r+1}$, we have from Properties \ref{p4} and \ref{p5}
\begin{align*}
    I(n) & \leq I\left(5^{2a_1} \cdot p_4^{2a_2} \cdot p_5^{2a_3} \dots \cdot p_{r+1}^{2a_{r-1}} \cdot \prod_{r\leq i\leq \omega(n)} p^{2a_i}_{ L+i-r}\right)\\
    &<\frac{5}{4}\cdot \prod_{4 \leq j\leq (r+1)} \frac{p_j}{p_j-1} \cdot \prod_{r \leq i\leq \omega(n)}\frac{p_{L+i-r}}{p_{L+i-r}-1}.
\end{align*}
Thanks to Remark 3.3 \cite{4} for confirming
\begin{align*}\label{31}
    I(n)&<\frac{5}{4}\cdot \prod_{4 \leq j\leq (r+1)} \frac{p_j}{p_j-1} \cdot \prod_{r \leq i\leq \omega(n)}\frac{L+i-r}{L+i-r-1}\\
    &=\frac{5}{4}\cdot \prod_{4 \leq j\leq (r+1)} \frac{p_j}{p_j-1} \cdot \frac{L+\omega(n)-r}{L-1}.\tag{1}
\end{align*}

We shall prove that for any $\omega(n)\in\mathbb{Z^+}$, the following inequality holds
\begin{align*}
   \frac{L+\omega(n)-r}{L-1}< \frac{\mathcal{A+B}}{\mathcal{A}}.
\end{align*}
Since $\frac{\mathcal{A}}{\mathcal{B}}\in \mathbb{Q}^{+}\setminus \mathbb{Z^{+}}$, we can write $\mathcal{A}=\mathcal{B}k+\mu$, where $k,\mu\in \mathbb{Z^+}$ and $\mu\in \{1,\dots, \mathcal{B}-1\}$, also note that $(\mathcal{B},\mu)=1$ as $(\mathcal{A,\mathcal{B}})=1$. Therefore,
\begin{align*}
\frac{L+\omega(n)-r}{L-1}=\frac{(k+1)\omega(n)-r+\lceil\frac{\mu\omega(n)}{\mathcal{B}}\rceil}{k\omega(n)-1+\lceil\frac{\mu \omega(n)}{\mathcal{B}}\rceil}.    
\end{align*}
We now consider the following two cases, where we essentially observe the behavior of $$\frac{(k+1)\omega(n)-r+\lceil\frac{\mu\omega(n)}{\mathcal{B}}\rceil}{k\omega(n)-1+\lceil\frac{\mu \omega(n)}{\mathcal{B}}\rceil},$$ based on the divisibility of $\omega(n)$ by $\mathcal{B}$.\\

\textbf{Case-1:}\\
If $\mathcal{B}$ does not divide $\omega(n)$ then $\frac{\mu \omega(n)}{\mathcal{B}}\not\in \mathbb{Z^+}$ as $(\mathcal{B},\mu)=1$ and so
\begin{align*}
  \frac{(k+1)\omega(n)-r+\lceil\frac{\mu\omega(n)}{\mathcal{B}}\rceil}{k\omega(n)-1+\lceil\frac{\mu \omega(n)}{\mathcal{B}}\rceil} &=\frac{(k+1)\omega(n)-r+1+\frac{\mu\omega(n)}{\mathcal{B}}-\{\frac{\mu\omega(n)}{\mathcal{B}}\}}{k\omega(n)+\frac{\mu \omega(n)}{\mathcal{B}}-\{\frac{\mu\omega(n)}{\mathcal{B}}\}}\\
  &=\frac{(\mathcal{B}k+\mathcal{B}+\mu)\omega(n)+(1-r)\mathcal{B}-\mathcal{B}\{\frac{\mu\omega(n)}{\mathcal{B}}\}}{(\mathcal{B}k+\mu)\omega(n)-\mathcal{B}\{\frac{\mu\omega(n)}{\mathcal{B}}\}}\\
  &=\frac{(\mathcal{A}+\mathcal{B})\omega(n)+(1-r)\mathcal{B}-\mathcal{B}\{\frac{\mu\omega(n)}{\mathcal{B}}\}}{\mathcal{A}\omega(n)-\mathcal{B}\{\frac{\mu\omega(n)}{\mathcal{B}}\}}.
    \end{align*}
Note that, for any positive integer $Q$ such that $\mathcal{B}$ does not divide $Q$, we can write $Q=\mathcal{B}q+v$, where $q\in \mathbb{Z}_{\geq 0}$ and $v\in \{1,2,\dots,\mathcal{B}-1\}$. Therefore, choosing $Q=\mu\omega(n)$, we get $\{\frac{\mu\omega(n)} {\mathcal{B}}\}\in \{\frac{1}{\mathcal{B}},\frac{2}{\mathcal{B}},\dots, \frac{\mathcal{B}-1}{\mathcal{B}}\}$ and thus,

\begin{align*}
    1\leq \mathcal{B}\{\frac{\mu\omega(n)}{\mathcal{B}}\}\leq \mathcal{B}-1
\end{align*}
that is,
\begin{align*}\label{3a}
    (\mathcal{A}+\mathcal{B})\omega(n)+(1-r)\mathcal{B}-\mathcal{B}\{\frac{\mu\omega(n)}{\mathcal{B}}\}\leq (\mathcal{A}+\mathcal{B})\omega(n)+(1-r)\mathcal{B}-1 \tag{2}
\end{align*}
    and
    \begin{align*}\label{3b}
        \mathcal{A}\omega(n)-\mathcal{B}+1\leq \mathcal{A}\omega(n)-\mathcal{B}\{\frac{\mu\omega(n)}{\mathcal{B}}\}. \tag{3}
    \end{align*}

Using inequalities (\ref{3a}) and (\ref{3b}), we get
\begin{align*}
   \frac{(\mathcal{A}+\mathcal{B})\omega(n)+(1-r)\mathcal{B}-\mathcal{B}\{\frac{\mu\omega(n)}{\mathcal{B}}\}}{\mathcal{A}\omega(n)-\mathcal{B}\{\frac{\mu\omega(n)}{\mathcal{B}}\}}\leq \frac{(\mathcal{A}+\mathcal{B})\omega(n)+(1-r)\mathcal{B}-1}{\mathcal{A}\omega(n)-\mathcal{B}+1}.
\end{align*}
Let us define $\phi: [1,\infty)	\rightarrow\mathbb{R}$ by $\phi(t)=\frac{(\mathcal{A}+\mathcal{B})t+(1-r)\mathcal{B}-1}{\mathcal{A}t-\mathcal{B}+1}$. Observe that, 
$$\phi'(t)=\frac{\mathcal{AB}(r-2)+2\mathcal{A}+\mathcal{B}-\mathcal{B}^2}{(\mathcal{A}t-\mathcal{B}+1)^2}$$
as $\mathcal{AB}(r-2)+2\mathcal{A}+\mathcal{B}>\mathcal{B}^2$, we have $\phi'(t)>0$ and thus $\phi$ is a strictly increasing function of $t$ in $[1,\infty)$. Since $\lim_{t\rightarrow \infty} \phi(t)=\frac{\mathcal{A+B}}{\mathcal{A}}$ we have $\phi(t)<\frac{\mathcal{A+B}}{\mathcal{A}}$ for all $t\in [1,\infty)$. In particular, setting $t=\omega(n)$ we get
\begin{align*}
\frac{(\mathcal{A}+\mathcal{B})\omega(n)+(1-r)\mathcal{B}-1}{\mathcal{A}\omega(n)-\mathcal{B}+1}<\frac{\mathcal{A+B}}{\mathcal{A}},
\end{align*}
that is, 
$$\frac{L+\omega(n)-r}{L-1}< \frac{\mathcal{A+B}}{\mathcal{A}}.$$
\textbf{Case-2:}\\
Let us consider the case, when $\mathcal{B}$ divides $ \omega(n)$, then $\frac{\mu \omega(n)}{\mathcal{B}}\in \mathbb{Z^+}$. Therefore,
\begin{align*}
 \frac{(k+1)\omega(n)-r+\lceil\frac{\mu\omega(n)}{\mathcal{B}}\rceil}{k\omega(n)-1+\lceil\frac{\mu \omega(n)}{\mathcal{B}}\rceil}=\frac{(k+1)\omega(n)-r+\frac{\mu\omega(n)}{\mathcal{B}}}{k\omega(n)-1+\frac{\mu \omega(n)}{\mathcal{B}}}&=\frac{(\mathcal{B}k+\mu+\mathcal{B})\omega(n)-\mathcal{B}r}{(\mathcal{B}k+\mu)\omega(n)-\mathcal{B}}\\
 &=\frac{(\mathcal{A}+\mathcal{B})\omega(n)-\mathcal{B}r}{\mathcal{A}\omega(n)-\mathcal{B}}.
 \end{align*}
Let us define $\tau: [1,\infty)	\rightarrow\mathbb{R}$ by $\tau(t)=\frac{(\mathcal{A}+\mathcal{B})t-\mathcal{B}r}{\mathcal{A}t-\mathcal{B}}$. Note that
$$\tau'(t)=\frac{\mathcal{AB}r-\mathcal{AB}-\mathcal{B}^2}{(\mathcal{A}t-\mathcal{B})^2}$$
as $r\geq 2$ and $\mathcal{A}>\mathcal{B}$, $\mathcal{AB}(r-1)-\mathcal{B}^2> \mathcal{AB}(r-1)-\mathcal{AB}=\mathcal{AB}(r-2)\geq 0$, it follows that $\tau'(t)>0$ and thus $\tau$ is a strictly increasing function of $t$ in $[1,\infty)$. Since $\lim_{t\rightarrow \infty} \tau(t)=\frac{\mathcal{A}+\mathcal{B}}{\mathcal{A}}$ we have $\tau(t)<\frac{\mathcal{A+B}}{\mathcal{A}}$ for all $t\in [1,\infty)$. In particular, setting $t=\omega(n)$ we get
\begin{align*}
\frac{(\mathcal{A}+\mathcal{B})\omega(n)-\mathcal{B}r}{\mathcal{A}\omega(n)-\mathcal{B}}<\frac{\mathcal{A+B}}{\mathcal{A}} 
\end{align*}
which immediately implies that 
\begin{align*}
   \frac{L+\omega(n)-r}{L-1}< \frac{\mathcal{A+B}}{\mathcal{A}}.
\end{align*}
Thus, for $\omega(n)\in \mathbb{Z^+}$, we have
\begin{align}\label{34}
  \frac{L+\omega(n)-r}{L-1}< \frac{\mathcal{A+B}}{\mathcal{A}}.\tag{4}
\end{align}
From (\ref{31}) and (\ref{34}), we get
\begin{align*}
    I(n)&<\frac{5}{4}\cdot \prod_{4 \leq j\leq (r+1)} \frac{p_j}{p_j-1} \cdot \frac{\mathcal{A+B}}{\mathcal{A}}\\
    &=\frac{5}{4}\cdot \prod_{4 \leq j\leq (r+1)} \frac{p_j}{p_j-1}\cdot(1+ \frac{\mathcal{B}}{\mathcal{A}})
\end{align*}
since we are given,
$$\frac{\mathcal{A}}{\mathcal{B}} > \frac{1}{\frac{36}{25}\cdot{\displaystyle \prod_{4\leq i\leq r+1} (1-\frac{1}{p_i})-1}},$$
we get
\begin{align*}
    I(n)&< \frac{5}{4}\cdot \prod_{4 \leq j\leq (r+1)} \frac{p_j}{p_j-1}\cdot(1+ \frac{\mathcal{B}}{\mathcal{A}})\\
    &<\frac{5}{4}\cdot \prod_{4 \leq j\leq (r+1)} \frac{p_j}{p_j-1}\cdot \frac{36}{25}\cdot \prod_{4\leq i\leq r+1} (1-\frac{1}{p_i})=\frac{9}{5}.
\end{align*}
Therefore, $n$ cannot be a friend of $10$ as $I(10)=\frac{9}{5}$. Hence, we must have $q_r < p_{L}$. This completes the proof. \qed

\section{Alternative Proof of Theorem \ref{lem1}}
For $q_2$, choose $\mathcal{A}=7$ and $\mathcal{B}=3$, then 
$\frac{\mathcal{A}}{\mathcal{B}}>\frac{25}{11}$, therefore applying Theorem \ref{thm 1.4} we get
$$q_2 < {\lceil\frac{7\omega(n)}{3}\rceil}\biggl\{\log\left({\lceil\frac{7\omega(n)}{3}\rceil}\right)+2\log\log\left({\lceil\frac{7\omega(n)}{3}\rceil}\right)\biggl\}.$$
For $q_3$, choose $\mathcal{A}=180$ and $\mathcal{B}=41$, then 
$\frac{\mathcal{A}}{\mathcal{B}}>\frac{175}{41}$, therefore applying Theorem \ref{thm 1.4} we get
$$q_3 < {\lceil\frac{180\omega(n)}{41}\rceil}\biggl\{\log\left({\lceil\frac{180\omega(n)}{41}\rceil}\right)+2\log\log\left({\lceil\frac{180\omega(n)}{41}\rceil}\right)\biggl\}.$$
For $q_4$, choose $\mathcal{A}=390$ and $\mathcal{B}=47$, then 
$\frac{\mathcal{A}}{\mathcal{B}}>\frac{385}{47}$, therefore applying Theorem \ref{thm 1.4} we get
$$q_4<  {\lceil\frac{390\omega(n)}{47}\rceil}\biggl\{\log\left({\lceil\frac{390\omega(n)}{47}\rceil}\right)+2\log\log\left({\lceil\frac{390\omega(n)}{47}\rceil}\right)\biggl\}.$$
This completes the proof.\qed\\
\section{Proof of Corollary \ref{coro1.4}}
For $q_3$, set $\mathcal{A}=427$ and $\mathcal{B}=100$, then 
$\frac{\mathcal{A}}{\mathcal{B}}>\frac{175}{41}$, therefore applying Theorem \ref{thm 1.4} we get
$$q_3 < {\lceil\frac{427\omega(n)}{100}\rceil}\biggl\{\log\left({\lceil\frac{427\omega(n)}{100}\rceil}\right)+2\log\log\left({\lceil\frac{427\omega(n)}{100}\rceil}\right)\biggl\}.$$
For $q_4$, set $\mathcal{A}=41$ and $\mathcal{B}=5$, then 
$\frac{\mathcal{A}}{\mathcal{B}}>\frac{385}{47}$, therefore applying Theorem \ref{thm 1.4} we get
$$q_4<  {\lceil\frac{41\omega(n)}{5}\rceil}\biggl\{\log\left({\lceil\frac{41\omega(n)}{5}\rceil}\right)+2\log\log\left({\lceil\frac{41\omega(n)}{5}\rceil}\right)\biggl\}.$$
This completes the proof.\qed\\
\section{Conclusion}
In order to find a friend $n$ of $10$ having exactly $m$ distinct prime divisors, we can apply Theorem \ref{thm 1.4} to get upper bounds for all $m$ distinct prime divisors of $n$, thereafter, we can consider all the possible combinations of prime divisors of $n$ one by one to check whether any of the possible combinations of prime divisors of $n$ can give abundancy index exactly equal to $\frac{9}{5}$. The solitariness status of $10, 14, 15, 20$ and many others is still unknown to us. Moreover, computer searches \cite{6} confirm that a friend of $10$ must be strictly larger than $10^{30}$.

\section*{Acknowledgment}
The author is grateful to the anonymous referees for the careful reading of the manuscript and for the valuable suggestions that greatly improved the clarity and presentation of the paper.

\end{document}